\documentclass[12pt]{article}
\usepackage{xcolor}
\usepackage{amsmath}
\usepackage{amssymb}
\usepackage{array}
\usepackage{float}

\newtheorem{proposition}{Proposition}
\newtheorem{theorem}{Theorem}
\newtheorem{corollary}{Corollary}

\begin{document}

\begin{center}{\Large
Symmetric measures of pseudorandomness for\\ binary sequences}
\end{center}

\begin{center} {\large Yixin Ren$^1$ and \large Arne Winterhof$^{2}$} \\

{\normalsize $^{1}$ Research Center for Number Theory and Its Applications,\\
School of Mathematics,} 
{\normalsize Northwest University
Xi’an 710127, China} \\
{\normalsize E-mail: 201920494@stumail.nwu.edu.cn}\\
{\normalsize $^{2}$ Johann Radon Institute for Computational and Applied Mathematics,} \\
{\normalsize Austrian Academy of Sciences, Altenbergestr. 69, A-4040 Linz, Austria} \\
{\normalsize E-mail: arne.winterhof@oeaw.ac.at}
 \end{center} 

\begin{abstract}
We compare ordinary and symmetric variants of two classical measures of pseudorandomness for binary sequences, the  $2$-adic complexity and the linear complexity.

In the periodic setting, we show that for binary periodic sequences constructed from the binary expansions of non-palindromic primes, the symmetric $2$-adic complexity can be strictly smaller than the ordinary $2$-adic complexity. We also give a direct proof (of the known result) that the linear complexity of a periodic binary sequence is invariant under reversal, and hence coincides with its symmetric version.

In the aperiodic setting, we provide explicit families of finite binary sequences for which both the 
$N$th symmetric 2-adic complexity and the $N$th symmetric linear complexity are substantially smaller than their ordinary counterparts. Furthermore, we show that the expected values of the 
$N$th rational complexity and of the 
$N$th exponential linear complexity exceed those of their symmetric analogues by at least a term of order of magnitude $N$. 
Thus, the effect of symmetrization is clearly visible on an exponential scale.

We also establish lower bounds for the expected values of the symmetric rational complexity, symmetric $2$-adic complexity, symmetric linear complexity, and symmetric exponential linear complexity.





\end{abstract}


\section{Introduction}
\subsection{Background}

Measures of pseudorandomness for binary sequences play a central role in the analysis and design of stream ciphers. Two of the most important complexity measures are linear complexity and $2$-adic complexity. Linear complexity is the length of the shortest linear feedback shift register (LFSR) that generates a given sequence. Its probabilistic theory was developed by Rueppel \cite{R86} and Niederreiter \cite{N88,N90}. In particular, the expected linear complexity of a random binary sequence of length $N$ is approximately $N/2$.

Klapper and Goresky \cite{GK12,KG97} introduced feedback with carry shift registers (FCSRs), whose output sequences are related to $2$-adic expansions of rational numbers. This led to the notion of $2$-adic complexity, which measures resistance to attacks based on FCSR synthesis. Later, Tian and Qi \cite{TQ10}, and more recently Chen and Winterhof \cite{CW25}, established corresponding expected-value results for finite binary sequences, showing that the expected 
$N$th 2-adic complexity is also of order $N/2$ up to a deviation of order of magnitude $\log(N)$.

In recent years, symmetric variants of these measures, obtained by minimizing over a sequence and its reversal, have been studied for certain special families of periodic sequences, such as generalized cyclotomic constructions, see for example \cite{EW25,ES22,EZ25}. However, a systematic comparison between ordinary and symmetric complexity measures, especially in the aperiodic setting, appears to be missing. The goal of this paper is to provide such a comparison for both the $2$-adic and the linear complexity, in periodic and aperiodic contexts.

 \subsection{Our results}
Let $\mathcal{S}=(s_n)_{n=0}^\infty$ be a binary sequence, that is, $s_n\in\{0,1\}$, $n=0,1,\ldots$
For $N=1,2,\ldots$ we put $\mathcal{S}_N=(s_0,s_1,\ldots,s_{N-1})$.

We study symmetric measures of pseudorandomness, first for periodic sequences and then for aperiodic sequences or finite parts of a periodic sequence, respectively.

In the periodic setting, our main results are:
\begin{itemize}
\item We derive sequences from non-palindromic primes with smaller symmetric $2$-adic complexity than $2$-adic complexity, see Theorem~\ref{prop1}.  
\item We prove directly from the recurrence definition of the linear complexity that the linear complexities of a periodic sequence and its reverse coincide. Consequently, the symmetric linear complexity coincides with the ordinary linear complexity in the periodic case, see Proposition~\ref{lin=sym}.
\end{itemize}

In the aperiodic setting,  we obtain: 
\begin{itemize}
\item families of binary sequences of substantially smaller symmetric $2$-adic complexity than $2$-adic complexity and substantially smaller symmetric linear complexity than the linear complexity;
\item an upper bound showing that the expected value of the $N$th rational complexity exceeds the expected value of the $N$th symmetric rational complexity by at least a summand of order of magnitude $N$, see Proposition~\ref{N-improve};
\item lower bounds for the expected values of the symmetric rational complexity, symmetric $2$-adic complexity, symmetric linear complexity, and symmetric exponential linear complexities, see Theorems~\ref{thm1} and \ref{thm3}.
\end{itemize}

\section{The periodic case}

First we consider sequences $\mathcal{S}$ of period $T$, that is, $s_{n+T}=s_n$, $n=0,1,\ldots$
and $\mathcal{S}$ is uniquely defined by the initial vector $\mathcal{S}_T$.

\subsection{$2$-adic complexity}

The {\em $2$-adic complexity} $\lambda(\mathcal{S})$ of $\mathcal{S}$ is 
$$\lambda(\mathcal{S})=\log_2\left(\frac{2^T-1}{\gcd(2^T-1,S_T(2))}\right),$$
where
$$S_T(x)=\sum_{n=0}^{T-1}s_nx^n.$$
 The $2$-adic complexity was introduced by Goresky and Klapper \cite{GK12,KG97} and is closely related
to the length of an
FCSR which generates the sequence.

The {\em reverse sequence} $\mathcal{S}^{rev}$ of $\mathcal{S}$ is the $T$-periodic sequence  defined by the initial vector 
$$\mathcal{S}_T^{rev}=(s_{T-1},s_{T-2},\ldots,s_0).$$
Obviously, we have
$$\lambda(\mathcal{S}^{rev})=\log_2\left(\frac{2^T-1}{\gcd(2^T-1,S_T^{rev}(2))}\right),$$
where 
$$S_T^{rev}(x)=\sum_{n=0}^{T-1}s_{T-1-n}x^n.$$
The {\em symmetric $2$-adic complexity} $\lambda^{sym}(\mathcal{S})$ is
$$\lambda^{sym}(\mathcal{S})=\min\{\lambda(\mathcal{S}),\lambda(\mathcal{S}^{rev})\}.$$

The symmetric $2$-adic complexity can be smaller than the $2$-adic complexity.

For example, for $T\ge 4$
consider the $T$-periodic sequence with 
$$\mathcal{S}_T=(1,1,0,1,0,0,\ldots,0),$$ that is, 
$$S_T(2)=\sum_{n=0}^{T-1}s_n2^n=11$$ 
and
$$S_T^{rev}(2)=\sum_{n=0}^{T-1}s_{T-1-n}2^n=2^{T-1}+2^{T-2}+2^{T-4}=13\cdot 2^{T-4}.$$
Now assume that $T\equiv 0\bmod 12$, that is, 
$$\gcd(2^T-1,13\cdot 2^{T-4})=\gcd(2^T-1,13)=13$$ 
and $T\not\equiv 0\bmod 10$, that is, $\gcd(2^T-1,11)=1$. Then we have
$$\lambda(\mathcal{S})=\lambda^{sym}(\mathcal{S})+\log_2(13)=\log_2(2^T-1).$$

This idea can be generalized as follows.
Two odd integers $p$ and $q$ are {\em reversible} (in base $2$) if 
\begin{equation}\label{pq}
p=\sum_{n=0}^{t-1} s_n2^n\quad \mbox{and}\quad q=\sum_{n=0}^{t-1} s_{t-1-n}2^n\quad \mbox{with } s_0=s_{t-1}=1
\end{equation}
and any $s_1,\ldots,s_{t-2}\in\{0,1\}$.
The integer $p$ is {\em palindromic} if $p=q$. 
We denote by ord$_p(2)$ the order of $2$ modulo $p$, that is, the smallest positive integer $o$ with $2^o\equiv 1\bmod p$.

\begin{theorem}\label{prop1}
    Let $p>2$ be a non-palindromic prime and write
    $$p=\sum_{n=0}^{t-1}s_n2^n, \quad s_0=s_{t-1}=1,~s_1,\ldots,s_{t-2}\in\{0,1\}.$$
    Let $q$ be the reverse of~$p$ defined by \eqref{pq} and $T$ be a positive integer with
    \begin{equation}\label{Tcond}
    T\equiv 0\bmod {\rm ord}_q(2)\quad \mbox{and}\quad T\not\equiv 0\bmod {\rm ord}_p(2).
    \end{equation}    
    Then the $T$-periodic sequence $\mathcal{S}$ with initial vector 
    $$\mathcal{S}_T=(s_0,s_1,\ldots,s_{t-1},0,0,\ldots,0)$$ 
    satisfies
    $$\lambda(\mathcal{S})=\lambda^{sym}(\mathcal{S})+\log_2(q)=\log_2(2^T-1).$$
\end{theorem}
Proof.
Condition \eqref{Tcond} implies
$$2^T-1\equiv 0\bmod q \quad\mbox{and}\quad 2^T-1\not\equiv 0\bmod p.$$
We get
$$\gcd(2^T-1,q)=q\quad \mbox{and}\quad \gcd(2^T-1,p)=1$$
since $p$ is a prime.
Since
$$S_T(2)=\sum_{n=0}^{t-1}s_n2^n=p$$
and
$$S_T^{rev}(2)=\sum_{n=T-t}^{T-1}s_{T-1-n}2^n=2^{T-t}q,$$
we get

$$\lambda(\mathcal{S})
   =\log_2\left(\frac{2^T-1}{\gcd(2^T-1,p)}\right)\\
   =\log_2\left(2^T-1\right),
$$
$$   \lambda(\mathcal{S}^{rev})
   =\log_2\left(\frac{2^T-1}{\gcd(2^T-1,2^{T-t}q)}\right)\\
   =\log_2\left(2^T-1\right)-\log_2(q)
$$
and the result follows.
\hfill $\Box$\\

Remarks. 1. Note that the number of palindromic primes $p$ with $2<p<2^t$ can be trivially upper bounded by $2^{t/2}$ since the bits $s_1,\ldots,s_{\lfloor (t-1)/2\rfloor}$ uniquely define $p$. This size is negligible compared to the number 
of all primes in this interval, that is, $\frac{2^t(1+o(1))}{t\log 2}$  
 by the prime number theorem. 
 Here, as usual, $o(1)$ stands for a function $f(t)$ with $\lim\limits_{t\rightarrow \infty}f(t)=0$.
 
2. If $q$ is also a prime, we may interchange the roles of $p$ and $q$ in Theorem~\ref{prop1}.
A list of orders of the smallest pairs of reversible primes $(p,q)$ with $p<q$ is given in Table~\ref{tab1}. 
Note that $(q,p)$ is another reversible prime pair, of course.
    \begin{table}[ht]
\centering
\begin{tabular}{c|c|c|c}
\hline
$p$ & $q$ & $\operatorname{ord}_p(2)$ & $\operatorname{ord}_q(2)$ \\
\hline
11 & 13 & 10 & 12 \\
23 & 29 & 11 & 28 \\
37 & 41 & 36 & 20 \\
43 & 53 & 14 & 52 \\
47 & 61 & 23 & 60 \\
67 & 97 & 66 & 48 \\
71 & 113 & 35 & 28 \\
83 & 101 & 82 & 100 \\
131 & 193 & 130 & 96 \\
151 & 233 & 15 & 29 \\
163 & 197 & 162 & 196 \\
167 & 229 & 83 & 76 \\
173 & 181 & 172 & 180 \\
199 & 227 & 99 & 226 \\
223 & 251 & 37& 50\\
\hline
\end{tabular}
\caption{Orders of reversible prime pairs $(p,q)$ with $p<q$}\label{tab1}
\end{table}

 3. Dartyge et al.\ \cite{DMRSS24} studied the distribution of reversible primes.
 More precisely, 
 we have
 $2^{t-1}< p,q<2^t$ since the first bit $s_0$ and the last bit $s_{t-1}$ are both $1$. For $t=1,2,\ldots$ let $\Theta(t)$ be the number of primes $p$
 with $2^{t-1}< p<2^t$ such that $q$ defined by \eqref{pq} is also a prime.
 By \cite[Corollary~1.2]{DMRSS24} we have
 $$\Theta(t)=O\left(\frac{2^t}{t^2}\right),$$
 where $f(t)=O(g(t))$ means $|f(t)|\le c g(t)$ for some constant $c>0$.
 In \cite[Section~8]{DMRSS24} the authors presented numerical investigations and a heuristic
argument that permitted them to formulate the conjecture
$$\Theta(t)=(3+o(1))\frac{2^{t-1}}{t^2}.$$
Recall that the number of palindromic~$p$ with $2^{t-1}< p<2^t$ is at most $2^{t/2}$ and thus negligible. 

4. Table~\ref{tab2} complements Table~\ref{tab1} by a list of reversible pairs $(p,q)$ with prime~$p$ and composite $q$.

 \begin{table}[ht]
\centering
\begin{tabular}{c|c|c|c}
\hline
$p$ & $q$ & $\operatorname{ord}_p(2)$ & $\operatorname{ord}_q(2)$ \\
\hline
19 & 25 & 18 & 20 \\
59 & 55 & 58 & 20 \\
79 & 121 & 13 & 110 \\
89 & 77 & 11 & 30 \\
103 & 115 & 17 & 44 \\
109 & 91 & 108 & 12 \\
137 & 145 & 8 & 28 \\
139 & 209 & 23 & 90 \\
149 & 169 & 25 & 156 \\
157 & 185 & 39 & 36 \\
179 & 205 & 60 & 20 \\
191 & 253 & 30 & 110 \\
211 & 203 & 70 & 84 \\
239 & 247 & 7 & 36 \\
\hline
\end{tabular}
\caption{Orders of reversible pairs $(p,q)$ with prime $p$ and composite $q$}\label{tab2}
\end{table}

5. Another family of sequences for which the symmetric $2$-adic complexity is smaller than the $2$-adic complexity consists of the reversal $\ell$-sequences with prime connection integer.
More precisely, any sequence $\mathcal{S}^{rev}$ with (least) period $T>1$ and connection integer
$q=\frac{2^T-1}{\gcd(S_T^{rev}(2),2^T-1)}$
satisfies $T={\rm ord}_q(2)$, see \cite[Corollary 4.2.5]{GK12}, and $\lambda(\mathcal{S}^{rev})=\log_2(q)$.
If $q$ is a prime such that $2$ is a primitive root modulo $q$, then $T=q-1$ is maximal
and $\mathcal{S}^{rev}$ is an {\em $\ell$-sequence}, see for example \cite[Definition~13.1.1]{GK12}.
By \cite[Theorem 16.2.5]{GK12} the reversal $\mathcal{S}=(\mathcal{S}^{rev})^{rev}$ of $\mathcal{S}^{rev}$ is not an $\ell$-sequence and must have larger $2$-adic complexity than $\mathcal{S}^{rev}$
since otherwise $r=\frac{2^T-1}{\gcd(S_T(2),2^T-1)}\le q$ would imply $q-1=T= {\rm ord}_r(2)\le r-1\le q-1$, that is, $r=q$ and $\mathcal{S}$ is an $\ell$-sequence, a contradiction.

For example, see \cite[Example 16.2.2]{GK12}, take the $\ell$-sequence $\mathcal{S}^{rev}$ of (least) period~$18$
with $S_{18}^{rev}(2)=220752$ and $\lambda(\mathcal{S}^{rev})=\log_2(19)$. Its reversal $\mathcal{S}$
is defined by $S_{18}(2)=10731$ and satisfies $\lambda(\mathcal{S})=\log_2(171)$.\\

6. It is clear that $\lambda(\mathcal{S})=\lambda(\mathcal{S}^{rev})$ for any $T$-periodic sequence $\mathcal{S}$ with palindromic $S_T(2)$, that is, $\mathcal{S}^{rev}=\mathcal{S}$. 
We also have $\lambda(\mathcal{S})=\lambda(\mathcal{S}^{rev})$ in some other cases, for example, if $S_T(2)=2^kq$ with a palindromic $q$ with $2^{T-k-1}<q<2^{T-k}$, that is,
${\cal S}_T=(\underbrace{00\ldots 0}_k1q_1 q_2\ldots q_2 q_1 1)$.
Then we have $S_T^{rev}(2)=q$ and $\gcd(2^T-1,S_T(2))=\gcd(2^T-1,S_T^{rev}(2))$.

Another example is given if $S_T(2)=2^k-1+2^kq$  with a palindromic $q<2^{T-k-1}$ and thus $S_T^{rev}(2)=q+2^{T-k}(2^k-1)$.

We get the same result again for any cyclic shift of these sequences, of course.

We will see later that although for these sequences there is no difference between the $2$-adic complexity and the symmetric $2$-adic complexity in the periodic case, in the aperiodic case we will get different values.

These examples shall show that sometimes the weakness of a sequence is already detected in the periodic case but sometimes we need the aperiodic measure defined below.

7. Note that in the case that $T$ is prime and $2^T-1$ is a Mersenne prime any non-constant $T$-periodic sequence has maximal $2$-adic complexity, which is another reason for the importance of studying aperiodic measures.

\subsection{Linear complexity}

The {\em linear complexity} $L(\mathcal{S})$ of $\mathcal{S}$ is
the smallest $L$ such that 
\begin{equation}\label{rec}\sum_{\ell=0}^L c_\ell s_{n+\ell}=0,\quad n=0,1,\ldots
\end{equation}
 By the $T$-periodicity of $\mathcal{S}$ we have 
$$L(\mathcal{S})\le T$$
and may restrict \eqref{rec} to the range $n=0,1,\ldots,2T-1-L\ge T-1$.

An equivalent definition is, see for example \cite[Lemma 8.2.1]{CDR04}, \begin{equation}\label{gcd}L(\mathcal{S})=\deg\left(\frac{x^T-1}{\gcd(x^T-1,S_T(x))}\right).
\end{equation}

We present a very short proof of \cite[Lemma~16.2.1]{GK12} directly based on the definition \eqref{rec} whereas the proof in \cite{GK12} uses \eqref{gcd}.

\begin{proposition}\label{lin=sym}
For any periodic sequence $\mathcal{S}$ we have $L(\mathcal{S}^{rev})=L(\mathcal{S})$.
\end{proposition}
Proof.
Let \eqref{rec}
be a linear recurrence of order $L=L(\mathcal{S})$ for the first $2T-L$ elements of the $T$-periodic sequence $\mathcal{S}$ with initial values $\mathcal{S}_T=(s_0,s_1,\ldots,s_{T-1})$. We may assume $L\le T$.
Put $$m=2T-1-n-L, \quad k=L-\ell \quad \mbox{and}\quad s_n'=s_{2T-1-n}$$ for $n=0,1,\ldots,2T-1$ and $k=0,1,\ldots,L$. Then we get
$$\sum_{k=0}^Lc_{L-k} s'_{m+k}=0,\quad m=0,1,\ldots,2T-L-1,$$
that is, $L(\mathcal{S})=L\ge L(\mathcal{S}^{rev})$. Since $(\mathcal{S}^{rev})^{rev}=\mathcal{S}$ we get the result. \hfill $\Box$\\

Consequently, there is no need to study the symmetric linear complexity 
$L^{sym}(\mathcal{S})=\min\{L(\mathcal{S}),L(\mathcal{S}^{rev})\}$
in the periodic case.

\section{The aperiodic case}

Let $\mathcal{S}=(s_n)_{n=0}^\infty$ be a binary sequence, that is, $s_n\in \{0,1\}$.
We write $\mathcal{S}_N=(s_0,s_1,\ldots,s_{N-1})$.
The {\em reverse} $\mathcal{S}_N^{rev}$ of $\mathcal{S}_N$ is
$$\mathcal{S}_N^{rev}=(s_{N-1},s_{N-2},\ldots,s_1,s_0).$$

\subsection{$2$-adic complexity}

The {\em $N$th rational complexity} $\Lambda(\mathcal{S}_N)$ of $\mathcal{S}$
is
$$\Lambda(\mathcal{S}_N)=\min\{\max\{q,|f|\}: q\sum_{n=0}^{N-1}s_n2^n\equiv f\bmod 2^N,~q>0 \mbox{ odd}\}$$
and the {\em $N$th $2$-adic complexity} $\lambda(\mathcal{S}_N)$ is
$$\lambda(\mathcal{S}_N)=\log_2(\Lambda(\mathcal{S}_N)),$$ 
see for example \cite{CW25,TQ10}.

 For $N\ge 2$, the {\em symmetric rational complexity is}
$$\Lambda^{sym}(\mathcal{S}_N)=\min\{\Lambda(\mathcal{S}_N),\Lambda(\mathcal{S}_N^{rev})\}$$
and the {\em symmetric $2$-adic complexity is}
$$\lambda^{sym}(\mathcal{S}_N)=\min\{\lambda(\mathcal{S}_N),\lambda(\mathcal{S}_N^{rev})\}.$$

We recall the results from \cite[Theorem 1]{CW25} on the expected values of the rational complexity and the $2$-adic complexity,
$$E_N^{rat}=\frac{1}{2^N}\sum_{\mathcal{S}_N\in \{0,1\}^N}\Lambda(\mathcal{S}_N)=2^{N/2+O(N/\log(N))},
$$
$$E_N^{2-adic}=\frac{1}{2^N}\sum_{\mathcal{S}_N\in \{0,1\}^N}\lambda(\mathcal{S}_N)=\frac{N}{2}+O(\log(N)).
$$

In this section we will estimate the symmetric analogs
$$E_N^{rat-sym}=\frac{1}{2^N}\sum_{\mathcal{S}_N\in \{0,1\}^N}\Lambda^{sym}(\mathcal{S}_N)$$
and
$$E_N^{2-adic-sym}=\frac{1}{2^N}\sum_{\mathcal{S}_N\in \{0,1\}^N}\lambda^{sym}(\mathcal{S}_N).$$

 The following examples show that the symmetric rational and the symmetric $2$-adic complexity can be substantially smaller than the rational and $2$-adic complexity, respectively.\\
 
Example $1$. For  $k\ge\frac{N}{2}$ define $\mathcal{S}_N=(s_0,s_1,\ldots,s_{N-1})$ by
$$s_n=0,\quad n=0,1,\ldots,k-1,\quad s_k=1,\quad (s_{k+1},s_{k+2},\ldots,s_{N-1})\in \{0,1\}^{N-k-1}.$$
By \cite[Lemma~4]{TQ10} 
we have 
$$\Lambda(\mathcal{S}_N)= 2^k,\quad k\ge \frac{N}{2}.$$
Obviously, we have (choosing $q=1$),
$$\Lambda(\mathcal{S}_N^{rev})\le  f$$
where
$$f=
\sum_{n=0}^{N-1}s_{N-1-n}2^n=\sum_{n=0}^{N-k-2}s_{N-1-n}2^n+2^{N-k-1}<2^{N-k}\le 2^{N/2}$$
and thus
$$\Lambda^{sym}(\mathcal{S}_N)=\Lambda(\mathcal{S}_N^{rev})<2^{N/2}\le \Lambda(\mathcal{S}_N).$$
For fixed $k$, $f$ runs through all integers between $2^{N-k-1}$ and $2^{N-k}-1$ whenever $(s_{k+1},s_{k+2},\ldots,s_{N-1})$ runs through $\{0,1\}^{N-k-1}$. Hence, $f$ runs through all integers between $1$ and $2^{\lfloor N/2\rfloor}-1$ when $k$ runs from $\lceil N/2\rceil$ to $N-1$.
So we get
$$E_N^{rat-sym}\le E_N^{rat}-M_1,$$
where 
\begin{eqnarray*}
M_1&=&2^{-N}\left(\sum_{k=\lceil N/2\rceil}^{N-1} 2^k \cdot 2^{N-k-1}
-\sum_{f=1}^{2^{\lfloor N/2\rfloor}-1} f\right)\\
&=&\frac{1}{2}\left(N-\left\lceil \frac{N}{2}\right\rceil\right)-\left(2^{\lfloor N/2\rfloor} -1\right)2^{\lfloor N/2\rfloor-N-1}=\frac{N-2}{4}+2^{-\lceil N/2\rceil-1}\\
&=&\frac{N}{4}+O(1).
\end{eqnarray*}

Example $2$. Now for $k\ge \frac{N+1}{2}$ consider
$$s_n=1,\quad n=0,1,\ldots,k-1,\quad s_k=0,\quad (s_{k+1},s_{k+2},\ldots,s_{N-1})\in \{0,1\}^{N-k-1}.$$
We have 
$$\Lambda(\mathcal{S}_N)\ge \Lambda(\mathcal{S}_k)\ge 2^{k-1}+1$$ 
since for any odd $q>0$,
$$q\sum_{n=0}^{N-1}s_n2^n\equiv q(2^k-1)\equiv 2^k-q\equiv f \bmod 2^{k+1}$$
and the minimum of $\max\{q,f\}$ is $2^{k-1}+1$.\\
We have 
$$\Lambda(\mathcal{S}_N^{rev})\le 2^{N-k}<\Lambda(\mathcal{S}_N)$$
since
$$\sum_{n=0}^{N-1}s_n2^{N-1-n}=2^N-2^{N-k}+\sum_{n=k+1}^{N-1}s_n2^{N-1-n}\equiv -2^{N-k}+S\equiv f\bmod 2^N$$
for some $S$ with $0\le S\le 2^{N-k-1}-1$.
Hence, for $q=1$ and $f=S-2^{N-k}$ we get
$$-f=|f|\le 2^{N-k}$$
and thus
$$\Lambda^{sym}(\mathcal{S}_N)\le 2^{N-k},\quad k\ge \frac{N+1}{2}.$$
Hence, we have $$E_N^{rat-sym}\le E_N^{rat}-M_1-M_2=E_N^{rat}-\frac{N-2}{4}-2^{-\lceil N/2\rceil -1}-M_2,$$
where
\begin{eqnarray*}
M_2&=&2^{-N} \sum_{k=\lceil (N+1)/2\rceil}^{N-1} 2^{N-k-1}\cdot (2^{k-1}+1-2^{N-k})\\
&=&\frac{N-\lceil (N+1)/2\rceil}{4}+2^{-\lceil (N+1)/2\rceil}-2^{-N} - \frac{1}{3}\left(2^{N-1-2\lceil (N+1)/2\rceil+2}-2^{-N+1}\right)\\
&=& \frac{N}{8}+O(1).
\end{eqnarray*}
We proved the following result.
\begin{proposition}\label{N-improve}
    $$E_N^{rat-sym}\le E_N^{rat}-\frac{3N}{8}+O(1).$$
\end{proposition}

Table \ref{tab:simple} compares the expected values 
$E_N^{rat}$ and $E_N^{rat-sym}$ and may lead to the conjecture that their difference is at least of  order of magnitude $2^{N/2}$.

    \begin{table}[h]  
    \centering
    \begin{tabular}{c c c c c c c c}  
     \hline  
       Length $N$  & 2 & 3 & 4 & 5 & 6 & 7 & 8 \\ 
        $E_N^{rat}-E_N^{rat-sym}$ & 0.25 & 0.5 & 0.813 & 1.344 & 1.922 & 2.813 & 3.992\\
       Length $N$   & 9 & 10 & 11 & 12 & 13& 14 & 15 \\ 
         $E_N^{rat}-E_N^{rat-sym}$  & 5.75 & 8.168 & 11.716 & 16.706 & 23.8 & 34.142 & 49.254\\
        Length $N$     & 16  & 17  & 18  & 19  & 20 & 21  \\
         $E_N^{rat}-E_N^{rat-sym}$  & 70.972 & 102.296 & 147.746 &213.142 & 307.238 & 441.608\\
        \hline  
    \end{tabular}
    \caption{Comparison of $E_N^{rat-sym}$ and $E_N^{rat}$}  
    \label{tab:simple}  
\end{table}

We will conclude this section with lower bounds on the expected values of the symmetric rational complexity and the symmetric $2$-adic complexity.


\begin{theorem}  \label{thm1}
For any function $r(N)$ with 
$$\lim_{N\rightarrow \infty}r(N)=\infty$$ 
we have $$E_N^{rat-sym}\geq 2^{N/2-r(N)+o(1)}.$$
Moreover, we have
$$E_N^{2-adic-sym}\geq \frac{N}{2}-\log(N)+o(1).$$
\end{theorem}
Proof.
By \cite[Lemma~1]{CW25}
with $W=2^{N/2-r(N)}$ the number of sequences $\mathcal{S}_N$ with $\Lambda(\mathcal{S}_N)\le W$ is $O(2^{N-2r(N)})$. Since each of these sequences appears at most twice in $\min\{\Lambda,\Lambda^{rev}\}$ we have for some $c>0$ at least $2^N(1-c4^{-r(N)})$ sequences with $\Lambda^{sym}\ge 2^{N/2-r(N)}$ and $\lambda^{sym}\ge \frac{N}{2}-r(N)$, respectively, thus 
$$E_N^{rat-sym}\ge (1-c4^{-r(N)})2^{N/2-r(N)}$$ 
and 
$$E_N^{2-adic-sym}\ge (1-c4^{-r(N)})\left(\frac{N}{2}-r(N)\right)$$
and the first result follows immediately.
For the second result we choose $r(N)=\log(N)$ and verify that $N4^{-\log(N)}=o(1)$.
\hfill $\Box$

\begin{corollary}
  We have
  $$E_N^{rat-sym}=2^{N/2+O(N/\log N)}$$
  and
  $$E_N^{2-adic-sym}=\frac{N}{2}+O(\log N).$$
\end{corollary}
Proof.
The result follows by combining Theorem~\ref{thm1} with the trivial estimates
$$E_N^{rat-sym}\le E_N^{rat} \quad \mbox{and}\quad E_N^{2-adic-sym}\le E_N^{2-adic}$$
and the upper bounds for $E_N^{rat}$ and $E_N^{2-adic}$ given in \cite[Theorem 1]{CW25}. \hfill $\Box$\\

Note that we know stronger lower bounds for $E_N^{rat}$ and $E_N^{2-adic}$ than the ones for their symmetric analogues given in Theorem~\ref{thm1}, see \cite[Theorem 2]{CW25},
$$E_N^{rat}\ge 2^{N/2-1}+\frac{N-5}{4}$$
and
$$E_N^{2-adic}\ge \frac{N}{2}-1.$$




\subsection{Linear complexity}

The {\em $N$th linear complexity} $L(\mathcal{S}_N)$ is the smallest $L$ such that
$$\sum_{\ell=0}^L c_\ell s_{n+\ell}=0,\quad n=0,1,\ldots,N-L-1,$$
and the {\em $N$th symmetric linear complexity}  is
$$L^{sym}(\mathcal{S}_N)=\min\{L(\mathcal{S}_N),L(\mathcal{S}_N^{rev})\}.$$
We recall the result on the expected value $E_N^{lin}$ of the $N$th linear complexity,
\begin{equation}\label{linexp}
E_N^{lin}=2^{-N}\sum_{\mathcal{S}_N\in \{0,1\}^N} L(\mathcal{S}_N)=\frac{N}{2}+O(1),
\end{equation}
see \cite{R86} or \cite[Section 18.1.1]{GK12}.

Now we also look into the expected values of $2^{L(\mathcal{S}_N)}$ and $2^{L^{sym}(\mathcal{S}_N)}$
to see the effect of the two families of sequences in the examples of the previous section.

We will use the well-known property that 
$$L(\mathcal{S}_{k+1})\in \{L(\mathcal{S}_k),k+1-L(\mathcal{S}_k)\},\quad k=0,1,\ldots$$
see for example \cite[Theorem 6.7.4]{J93}.

Example 3.
We consider again the sequence of Example 1 in the previous section, that is,
$$(s_0,s_1,\ldots,s_{k-1},s_k)=(0,0,\ldots,0,1)\quad \mbox{and} \quad (s_{k+1},\ldots,s_{N-1})\in \{0,1\}^{N-k-1}.$$ 
We get
$L(\mathcal{S}_k)=0$ and thus $L(\mathcal{S}_N)\ge L(\mathcal{S}_{k+1})=k+1-L(\mathcal{S}_k)=k+1$.\\
 Since $s^{rev}_{n+N-k}=s_{k-n-1}=0$ for $n=0,1,\ldots,k-1$ we get
$L(\mathcal{S}^{rev}_N)\le N-k$.
The contribution of these sequences to the difference of the expected values of $2^{L(\mathcal{S}_N)}$ and $2^{L^{sym}(\mathcal{S}_N)}$ is at least
\begin{equation}\label{el1}
K_1=2^{-N}\sum_{k=\lceil N/2\rceil}^{N-1}2^{N-1-k}(2^{k+1}-2^{N-k})=\frac{N}{2}+O(1).
\end{equation}

Example 4. We consider again the sequence of Example 2 in the previous section, that is, 
$$(s_0,s_1,\ldots,s_{k-1},s_k)=(1,1,\ldots,1,0)\quad \mbox{and} \quad (s_{k+1},\ldots,s_{N-1})\in \{0,1\}^{N-k-1}.$$
On the one hand, we get
$L(\mathcal{S}_k)=1$ since $s_{n+1}=s_n$ for $n=0,\ldots,k-2$
and 
$$L(\mathcal{S}_N)\ge L(\mathcal{S}_{k+1})=k+1-L(\mathcal{S}_k)=k.$$
On the other hand, we have 
$$L(\mathcal{S}_N^{rev})\le N-k+1$$ 
since $s^{rev}_{n+N-k+1}=s_{k-n-2}=s_{k-n-1}=s^{rev}_{n+N-k}=1$ for $n=0,\ldots,k-2$.
The difference between the expected value of $2^L$ and $2^{L^{rev}}$ is at least
\begin{equation}\label{el2}
K_2=2^{-N}\sum_{k=\lceil N/2\rceil}^{N-1}2^{N-k-1}(2^k-2^{N-k+1})=
\frac{N}{4}+O(1).
\end{equation}

For the expected values
$$E_N^{lin-exp}=2^{-N}\sum_{\mathcal{S}_N\in \{0,1\}^N} 2^{L(\mathcal{S}_N)}$$
and
$$E_N^{lin-exp-sym}=2^{-N}\sum_{\mathcal{S}_N\in \{0,1\}^N} 2^{L^{sym}(\mathcal{S}_N)}$$
we get the following result by combining \eqref{el1} and \eqref{el2}.
\begin{proposition}\label{prop3}
    $$E_N^{lin-exp-sym}\le E_N^{lin-exp}-\frac{3N}{4}+O(1).$$
\end{proposition}
Comparison of $E_N^{lin-exp-sym}$ and $E_N^{lin-exp}$ in Table~\ref{tab:last} may lead to the conjecture that their difference is of order of magnitude $2^{N/2}$.

\begin{table}[h]  
    \centering
    \begin{tabular}{c c c c c c c }  
     \hline  
       Length $N$  & 2 & 3 & 4 & 5 & 6 & 7 \\ 
        $E_N^{lin-exp}-E_N^{lin-exp-sym}$    & 0.5  & 0.75  & 1.625  &   2.063 & 3.656   & 4.453  \\
       Length $N$  & 8  & 9 & 10 & 11 & 12 & 13 \\ 
        $E_N^{lin-exp}-E_N^{lin-exp-sym}$   & 7.539   & 9.082  & 15.197  & 18.255  &   30.456 & 36.5559   \\
        Length $N$  & 14 & 15   & 16  & 17  & 18  & 19  \\
          $E_N^{lin-exp}-E_N^{lin-exp-sym}$   & 60.942  & 73.135   & 121.9  & 146.282  & 243.807  &   292.57  \\
          Length $N$  & 20 & 21 & 22 \\
          $E_N^{lin-exp}-E_N^{lin-exp-sym}$    & 487.618   & 585.142  & 975.237  \\
        \hline  
    \end{tabular}
    \caption{Comparison of $E_N^{lin-exp}$ and $E_N^{lin-exp-sym}$}  
    \label{tab:last} 
\end{table}


Now we prove lower bounds on the expected values $E_N^{lin-exp-sym}$ and~$E_N^{lin-sym}$.
\begin{theorem}\label{thm3}
For any function $r$ with 
$$\lim_{N\rightarrow \infty} r(N)=\infty$$
we have
$$E_N^{lin-exp-sym}\ge 2^{N/2-r(N)+o(1)}.$$
Moreover, we have
$$E_N^{lin-sym}\ge\frac{N}{2}-\log(N)+o(1).$$
\end{theorem}
Proof.
First we recall that the number $A_N(L)$ of sequences $\mathcal{S}_N\in \{0,1\}^N$
with $L(\mathcal{S}_N)=L$ is
$$
A_N(L) = \begin{cases} 
1 & \text{if } L(S_N) = 0, \\ 
2^{\min(2L-1, 2N - 2L)} & \text{if } 1 \leq L \leq N,
\end{cases}
$$
see for example \cite{GK12,R86}.


For $W\le \frac{N}{2}$ the number $M(W)$ of ${\cal S}_N$ with $L({\cal S}_N)\le W$ is
$$M(W)=\sum_{L=0}^WA_N(L)=1+\sum_{L=1}^W 2^{2L-1}=\frac{4^{W+1}+2}{6}=O(4^W).$$
With 
$$W=\frac{N}{2}-r(N)$$ we get
$$\frac{M(W)}{2^N}=O(4^{-r(N)}).$$
Each sequence $\mathcal{S}_N$ of small linear complexity $L$ appears at most twice as either $\mathcal{S}_N$ or $\mathcal{S}_N^{rev}$ in $L^{sym}(\mathcal{S}_N)=\min\{L(\mathcal{S}_N),L(\mathcal{S}_N^{rev})\}=L$
and thus as in the proof of Theorem~\ref{thm1} with some constant $c>0$,
\begin{eqnarray*}
E_N^{lin-exp-sym}&\ge& \left(1-\frac{2M(W)}{2^N}\right)2^{\frac{N}{2}-r(N)}\\
&\ge&\left(1-c4^{-r(N)}\right)2^{\frac{N}{2}-r(N)}
\end{eqnarray*}
and
\begin{eqnarray*}
E_N^{lin-sym}
&\ge& \left(1-\frac{2M(W)}{2^N}\right)\left(\frac{N}{2}-r(N)\right)\\
&\ge &\left(1-c4^{-r(N)}\right)\left(\frac{N}{2}-r(N)\right)
\end{eqnarray*}
and the results follow as in the proof of Theorem~\ref{thm1} with the
choice $r(N)=\log(N)$ for the second result. \hfill $\Box$


We complete our results with a formula for the expected value $E_N^{lin-exp}$. 
\begin{proposition}\label{lastprop}
We have
$$E_N^{lin-exp}= c2^{\frac{N}{2}}+O(1),$$
where 
\begin{equation}\label{cconst}
c=\left\{\begin{array}{cc} 11/7=1.57\ldots, & N \mbox{ even},\\
8\sqrt{2}/7=1.61\ldots, & N \mbox{ odd}.\end{array}\right.
\end{equation}
\end{proposition}
Proof.
%
We have
\begin{eqnarray*}
    E_N^{lin-exp} &=& 2^{-N}\sum_{L=1}^N 2^{\min\{2L-1,2N-2L\}}2^L\\
    &=&\frac{1}{2^N}\left(\sum_{L=1}^{\lfloor N/2\rfloor}2^{3L-1}+\sum_{L=\lfloor N/2\rfloor+1}^N2^{2N-L}\right)\\
    &=& \frac{1}{2^N}\left(\frac{4}{7}(2^{3\lfloor N/2\rfloor}-1)
    +2^{2N-\lfloor N/2\rfloor}-2^N\right)\\
    &=& c2^{N/2}-1-\frac{4}{7}2^{-N}\\
    &=& c2^{N/2}+O(1)
\end{eqnarray*}
and the result follows. \hfill $\Box$
%

Propositions \ref{prop3} and \ref{lastprop}, \eqref{linexp} and the trivial bound
$$E_N^{lin-sym}\le E_N^{lin}$$
provide the following result.
\begin{corollary}
    With the constant $c$ defined by \eqref{cconst} we have
    $$E_N^{lin-exp-sym}\le c2^{N/2}-\frac{3N}{4}+O(1)$$
    and
    $$E_N^{lin-sym}\le \frac{N}{2}+O(1).$$
\end{corollary}

\section*{Acknowledgments}
The work was written during a pleasant visit of Y. Ren to Linz. She wishes to thank the Chinese Scholarship Council for financial support and RICAM of the Austrian Academy of Sciences for hospitality.

The research of the second author was funded by the Austrian Science Fund (FWF) [10.55776/PAT4719224].

We'd like to thank one of the referees of \cite{CW25} for pointing to the problem of estimating expected values of symmetric measures.

We'd like to thank Zhixiong Chen for some useful comments.

\end{document}